\documentclass[12pt]{article}
\usepackage{pslatex}
\usepackage{amsmath,amsthm,amssymb,amsfonts, stmaryrd}
\usepackage{fancybox}
\usepackage{xcolor, graphicx}
\newtheorem{theorem}{Theorem}[section]

\newtheorem{definition}[theorem]{Definition}

\newtheorem{problem}[theorem]{Problem}
\newtheorem{proposition}[theorem]{Proposition}
\newtheorem{remark}[theorem]{Remark}

\newcommand{\Z}{\mathbb Z}

\newcommand{\N}{\mathbb N}
\newcommand{\PN}{\mathbb P}

\newcommand{\T}{\mathbb T}

\newcommand{\f}{\textit{fact}}

\date{}

\begin{document}

\title{Are monochromatic Pythagorean triples unavoidable under morphic colorings ?}
\author{S. Eliahou, J. Fromentin, V. Marion-Poty, D. Robilliard}

\maketitle

\begin{abstract} A Pythagorean triple is a triple of positive integers $a,b,c \in {\mathbb N}_+$ satisfying $a^2+b^2=c^2$. Is it true that, for any finite coloring of ${\mathbb N}_+$, at least one Pythagorean triple must be monochromatic? In other words, is the Diophantine equation $X^2+Y^2=Z^2$ \emph{regular}? This problem, recently solved for 2-colorings by massive SAT computations [Heule et al., 2016], remains widely open for $k$-colorings with $k \ge 3$. In this paper, we introduce \emph{morphic colorings} of ${\mathbb N}_+$, which are special colorings  in finite groups with partly multiplicative properties. We show that, for many morphic colorings in $2$ and $3$ colors, monochromatic Pythagorean triples are unavoidable in rather small integer intervals.
\end{abstract}

\medskip

\section{Introduction}

A triple $(a,b,c)$ of positive integers is a \emph{Pythagorean triple} if it satisfies $a^2+b^2=c^2$, as $(3,4,5)$ for instance. Is it true that, for any finite coloring of the set of positive integers, monochromatic Pythagorean triples are unavoidable?  While this typical Ramsey-type question has been open for decades \cite{EG:81}, there is no consensual conjecture as to whether the answer should be positive or not \cite{G:08}. Yet in the simplest case of two colors, the problem has just been solved in the affirmative with the support of massive SAT computations, as follows \cite{HKM:15}.

\begin{theorem}\label{heule} For any $2$-coloring of the integer interval $I=[1,7825]$, there is a monochromatic Pythagorean triple in $I$. Moreover, $7825$ is minimal with respect to this property.
\end{theorem}

Prior to this achievement, J. Cooper and R. Overstreet had obtained, already with SAT computations, an exotic $2$-coloring of the integer interval $[1,7664]$ avoiding monochromatic Pythagorean triples~\cite{CO:15}. See also~\cite{CFHW:14} for a related earlier work. The general problem remains widely open for $k$-colorings with $k \ge 3$. For background on Ramsey theory, see \cite{GRS:90}.

In this paper, we tackle the problem by focusing on a restricted class of  colorings that we call \textit{morphic}\footnote{We thank the anonymous referee for suggesting this term.}\textit{colorings}, which are partly multiplicative maps in some group and depend on the choice of a few prime numbers. For all the morphic colorings in $2$ and $3$ colors considered here, monochromatic Pythagorean triples turn out to be unavoidable, as expected in the case of $2$ colors of course, but much sooner so than the general threshold of $7825$ in Theorem~\ref{heule}.

The present results lead us to conjecture that the answer to the question in the title is, in fact, positive.

%
%
%
%
%

\section{Notation and background}

We shall denote by $\N$ the set of nonnegative integers, by $\N_+ = \{n \in \N \mid n \ge 1\}$ the subset of positive integers, and by $\PN=\{2, 3, 5, 7, 11, 13, \dots \}$ the subset of prime numbers. Given positive integers $a \le b$, we shall denote the integer interval they bound by 
$$
[a,b] \ = \ \{c \in \N \mid a \le c \le b\}.
$$

\begin{definition}
A \emph{Pythagorean triple} is a triple $(a,b,c)$ of positive integers satisfying $a^2+b^2=c^2$.  Such a triple is said to be \emph{primitive} if it satisfies $\gcd(a,b,c)=1$.
\end{definition}

Obviously, since the equation $X^2+Y^2=Z^2$ is homogeneous, every Pythagorean triple is a scalar multiple of a primitive one. 

The parametrization of primitive Pythagorean triples is well known. Indeed, \textit{every} primitive Pythagorean triple is of the form
$$
(m^2-n^2,2mn,m^2+n^2),
$$
where $m,n$ are coprime positive integers such that $m - n$ is positive and odd.

Following Rado \cite{Rad:33}, a Diophantine equation $f(X_1,\dots,X_n) = 0$ is said to be \textit{partition-regular}, or \textit{regular} for short, if for every finite coloring of $\N_+$, there is a monochromatic solution $(x_1,\dots,x_n) \in \N_+^n$ to it. More specifically, for given $k \in \N_+$, the equation is said to be \textit{$k$-regular} if, for every $k$-coloring of $\N_+$, there is a monochromatic solution to it. Note that regularity is equivalent to $k$-regularity for all $k \in \N_+$, and that $k$-regular implies $(k-1)$-regular if $k \ge 2$.

With this terminology, the question under study here is to determine whether the Diophantine equation $X^2+Y^2-Z^2=0$ is regular or not. And, if not, to determine the largest $k \ge 2$ for which this equation is $k$-regular. This problem is open since several decades. The only positive result about it so far is Theorem~\ref{heule}, which states that this equation is $2$-regular. It has been achieved by massive computations with a SAT solver in 35000 hours and validated in 16000 more hours with a 200TB certificate in the DRAT format \cite{HKM:15}.

The interest of our present approach with morphic colorings is three-fold. First, it allows us to test some $3$-colorings. Second, for all the $2$-colorings considered here, monochromatic Pythagorean triples turn out to be unavoidable much sooner, and at a much lower computational cost, than in Theorem~\ref{heule}. And third, as in Tao's recent solution of the Erd\H{o}s discrepancy conjecture~\cite{Tao:16}, partially multiplicative functions, like our morphic colorings, seem to be a good testbed for the problem under study.\footnote{Thanks are due to the anonymous referee for this observation.}

\section{Standard and partial morphisms }
\label{sect3}

In this paper, we introduce colorings of $\N_+$ by a finite additive group $G$ and satisfying special algebraic properties. Below, we shall mainly focus on the groups $\Z/2\Z$ and $\Z/3\Z$.

To start with, we consider monoid morphisms in the usual sense, i.e. maps $f \colon \N_+ \to G$ satisfying
$$
f(xy) = f(x)+f(y)
$$
for all $x,y \in \N_+$. Note that such a morphism is \textit{completely and freely determined by its values $\{f(p)\}_{p \in \PN}$ on the prime numbers}.

\smallskip

These morphisms $f$ are particularly interesting in the present context, since if a primitive Pythagorean triple $(a,b,c)$ fails to be monochromatic under $f$, then the same holds for all its scalar multiples $(ad,bd,cd)$ with $d \in \N_+$. Indeed, if $f(x)\not=f(y)$, then $f(xd) \not= f(yd)$ for all $d \ge 1$. This follows from the property $f(zd) = f(z)+f(d)$ for all $z \in \N_+$ and the fact that the values lie in a group.

\smallskip

We shall prove in a subsequent section that, for any morphism $f \colon \N_+ \to G$ where  $G$ is either $\Z/2\Z$ or $\Z/3\Z$, monochromatic primitive Pythagorean triples are unavoidable. Of course, for $G=\Z/2\Z$, that result follows from Theorem~\ref{heule}. However, its computer-aided proof is computationally much lighter; moreover, the unavoidability threshold turns out to be $533$ only in that case, as compared to $7825$ in the general case.

These positive results lead us to somewhat relax the constraints of morphisms if we seek to observe new phenomena, if any. Yet some structure on the considered colorings is needed, so as to have a more manageable function space size. This prompts us to consider maps $f \colon \N_+ \to G$ satisfying weaker conditions than morphisms, and which we now define.

\medskip

First, for any positive integer $n$, we denote by supp($n$) the set of prime factors of $n$. For instance, supp$(60)=\{2,3,5\}$.
 
\begin{definition} Let $(G,+)$ be an abelian group. Let $\PN_{0} \subseteq \PN$ be a given subset of the prime numbers. We say that a map
$$
f \colon \N_+ \to G
$$
is a \emph{$\PN_0$-partial morphism}, or a \emph{$\PN_0$-morphism} for short,  if the following properties hold. For any $n \in \N_+$, let $n_0 \in \N_+$ be the largest factor of $n$ such that supp$(n_0) \subseteq \PN_{0}$ and let $n_1=n/n_0$, so that $n=n_0n_1$. Then
\begin{itemize}
\item $f(n) = f(n_0) + f(n_1)$;
\item $f(n_1)=f(a)+f(b)$ for any coprime integers $a,b$ satisfying $n_1=a b$.
\end{itemize}
\end{definition} 
An equivalent way of expressing this notion is as follows. For any $n \in \N_+$, consider its unique prime factorization
$$
n = \prod_{p \in \PN} p^{\nu_p(n)}
$$
where $\nu_p(n) \in \N$ for all $p$. Then, the map $f \colon \N_+ \to G$ is a $\PN_0$-morphism if for all $n \in \N_+$, we have
$$
f(n) = f\big(\prod_{p \in \PN_0} p^{\nu_p(n)}\big) + \sum_{p \notin \PN_0} f\big(p^{\nu_p(n)}\big).
$$
Thus, a $\PN_0$-morphism is entirely and freely determined by its values on the set of positive integers
$$
S(\PN_0) = \big\{n_0 \in \N_+ \mid \textrm{supp}(n_0) \subseteq \PN_0\big\} \, \bigsqcup \, \big\{p^\nu \mid p \in \PN \setminus \PN_0, \nu \in \N_+\big\}.
$$
For instance, any $\{2,3\}$-morphism is freely  determined by its values on the integers of the form $2^a3^b$ and $p^c$ with $p \in \PN$, $p \ge 5$, where $a,b,c\in \N$ and $a+b \ge 1$, $c \ge 1$.

\begin{remark} Here are a few easy observations about $\PN_0$-morphisms $f \colon \N_+ \to G$.
\begin{itemize}
\item If $\PN_0 = \emptyset$, then $f$ is characterized by the property $f(xy)=f(x)+f(y)$ for all coprime positive integers $x,y$. In particular, classical monoid morphisms are $\emptyset$-morphisms in the present sense.
\item On the other end of the spectrum, if $\PN_0=\PN$, then $f$ is just a set-theoretical map without any special property or structure.
\item More generally, if $\PN_0 \subseteq \PN_1 \subseteq \PN$, then any $\PN_0$-morphism is also a $\PN_1$-morphism.
\end{itemize}
\end{remark}

\section{Colorings by standard morphisms}

The interest of using morphisms $f \colon \N_+ \to G$ as coloring functions is that such a coloring admits a monochromatic Pythagorean triple if and only if it admits a monochromatic \emph{primitive} Pythagorean triple. This is why we only need consider primitive Pythagorean triples in this section.

\begin{proposition}\label{morphism m=2} For any morphism $f \colon \N_+ \to \Z/2\Z$, monochromatic primitive Pythagorean triples are unavoidable. More precisely, such monochromatic triples are already unavoidable in the integer interval $[1,533]$. And finally, 533 is minimal with respect to this property.
\end{proposition}
The proof below relies on some computer assistance but, with patience, everything can be checked by hand. 
\begin{proof} Let $f \colon \N_+ \to \Z/2\Z$ be any morphism. Then $f$ is determined by its values on the prime numbers via the formula
$$
f(n)=f\big(\prod_{p \in \PN} p^{\nu_p(n)}\big) = \sum_{p \in \PN} \nu_p(n) f(p)
$$
for any $n \in \N_+$. Plainly, the only primes $p$ which actually contribute to the value of $f(n)$ are those for which $\nu_p(n)$ is odd. For example, we have $f(12)=f(3)$. 

For $n \in \N_+$, let us denote by oddsupp($n$) the \textit{odd support} of $n$, i.e. the set of primes $p$ for which $\nu_p(n)$ is odd. Thus, the above formula for $f(n)$ reduces to
$$
f(n)= \sum_{p} f(p),
$$
where $p$ runs through oddsupp$(n)$ only.

\smallskip

We shall restrict our attention to the 13 first primes, denoted $p_1,\dots,p_{13}$ in increasing order, and shall denote their set by $\PN_{13}$. Thus $\PN_{13}=\{2,3, \dots, 37,41\}$. Further, let us set
$$
\N_{|\PN_{13}} = \{n \in \N_+ \mid \textrm{oddsupp}(n) \subseteq \PN_{13} \}.
$$
For instance, the first few positive integers \textit{not} in $\N_{|\PN_{13}}$, besides the primes $p \ge 43$, are 86, 94, 106, 118, 122, etc.

\medskip

By the above formula, the value of $f(n)$ for any $n \in \N_{|\PN_{13}}$ is entirely determined by the length 13 binary vector
$$
v(f)=\big(f(p_1),\dots,f(p_{13})\big) \in (\Z/2\Z)^{13}.
$$
Consider now the set $\mathcal{T}$ of all primitive Pythagorean triples in the integer interval $[1,532]$. There are 84 of them, the lexicographically last one being $\{279, 440, 521\}$. Among them, we shall distinguish the subset $\mathcal{T}_{13}$ defined as
$$
\mathcal{T}_{13} = \big\{(a,b,c) \in \mathcal{T} \mid a,b,c \in \N_{|\PN_{13}} \big\},
$$
i.e. those triples in $\mathcal{T}$ whose three elements have odd support in $\PN_{13}$.
One finds that $|\mathcal{T}_{13}|=32$. For definiteness, here is this set:

\bigskip
\noindent
{\small
$\mathcal{T}_{13}=\big\{\{3,4,5\}$, $\{5,12,13\}$, $\{8,15,17\}$, $\{7,24,25\}$, $\{20,21,29\}$, $\{12,35,37\}$, $\{9,40,41\}$, $\{33,56,65\}$, $\{16,63,65\}$, $\{13,84,85\}$, $\{36,77,85\}$, $\{44,117,125\}$, $\{17,144,145\}$, $\{24,143,145\}$, $\{119,120,169\}$, $\{57,176,185\}$, $\{104,153,185\}$, $\{133,156,205\}$, $\{84,187,205\}$, $\{21,220,221\}$, $\{140,171,221\}$, $\{161,240,289\}$, $\{204,253,325\}$, $\{36,323,325\}$, $\{135,352,377\}$, $\{152,345,377\}$, $\{87,416,425\}$, $\{297,304,425\}$, $\{31,480,481\}$, $\{319,360,481\}$, $\{155,468,493\}$, $\{132,475,493\}\big\}$.
}

Perhaps surprisingly, it turns out that there are exactly two \textit{avoiding morphisms} 
$$
f_1,f_2 \colon \N_{|\PN_{13}} \to \Z/2\Z
$$
for which no $(a,b,c) \in \mathcal{T}_{13}$ is monochromatic. They are determined by the length 13 binary vectors $v(f_1)=w_1$, $v(f_2)=w_2$, where 
\begin{eqnarray*}
w_1 & = & 0101111101001,\\
w_2 & = & 0101111111001.
\end{eqnarray*}
Note that $w_1,w_2$ only differ at the 9th digit.

Now, the 85th primitive Pythagorean triple is $(308, 435, 533)$. As it happens, that triple is mapped to $(0,0,0)$ by both $f_1$ and $f_2$. Indeed, the prime factorizations of 308, 435 and 533 only involve the primes
$$
\begin{matrix}
p_{1}=2, & p_{2}=3, & p_{3}=5, & p_{4}=7,\\
p_{5}=11,& p_{6}=13, & p_{10}=29, & p_{13}=41,
\end{matrix}
$$
and are the following: $308=p_1^2p_4p_5$, $435=p_2p_3p_{10}$, $533=p_6p_{13}$. Hence, for $f = f_1$ or $f_2$, we have
$$
\begin{array}{lllllll}
f(308) & =  & f(p_{4})+f(p_{5}) & = & 1 + 1 & = & 0, \\
f(435) & =  & f(p_{2})+ f(p_{3})+f(p_{10})   & =  & 1 + 0 + 1 & =  & 0, \\
f(533) & =  & f(p_{6})+f(p_{13})  & =  & 1 + 1 & =  & 0.
\end{array}
$$
We conclude, as claimed, that for every morphism $g \colon \N_+ \to \Z/2\Z$, there must be a primitive Pythagorean triple in $[1,533]$ which is monochromatic under $g$. 

The fact that 533 is minimal with respect to this property is witnessed by the existence of many morphisms $f$ under which none of the 84 primitive Pythagorean triples in $[1,532]$ is monochromatic. The values on $p_1,\dots,p_{13}$ of these avoiding morphisms $f$ must of course be specified by either $w_1$ or $w_2$, but there are further restrictions. Indeed, their values on all primes turn out to be constrained as follows: either
$$
f(p_i) = \left\{
\begin{array}{ll}
0 & \textrm{for } i = 1,3,9,11,12,18,21,30,57,74,80,89, \\
1 & \textrm{for } i = 2,4,5,6,7,8,10,13,16,24,26,55,65,
\end{array}
\right.
$$
or
$$
f(p_i) = \left\{
\begin{array}{ll}
0 & \textrm{for } i = 1,3,11,12,18,21,25,30,59,74,89, \\
1 & \textrm{for } i = 2,4,5,6,7,8,9,24,26,55,65,70,
\end{array}
\right.
$$
with complete freedom on all other primes. Note that the 13 first bits of the first and of the second type, i.e. the $f(p_i)$'s for $i \le 13$, make up $w_1$ and $w_2$, respectively.
\end{proof}

\smallskip

An analogous result holds for morphic 3-colorings, established by an exhaustive computer search.

\begin{proposition}\label{morphism m=3} For any morphism $f \colon \N_+ \to \Z/3\Z$, monochromatic primitive Pythagorean triples are unavoidable. More precisely, at least one such triple in the integer interval $[1,4633]$ is monochromatic under $f$. And $4633$ is minimal with respect to that property.
\end{proposition}

Here is one particular morphism $f \colon \N_+ \to \Z/3\Z$ for which no Pythagorean triple in the interval $[1,4632]$ is monochromatic; it suffices to specify which primes in that interval are colored 1 or 2, the rest being colored 0. Denoting by $p_i$ the $i$th prime for $i \ge 1$, so that $p_1=2$, $p_2=3$, $p_3=5$ and so on, we set:
$$
f(p_i) \ = \ 
\left\{
\begin{array}{rl}
1 & \textrm{ if } i \in A, \\
2 & \textrm{ if } \ i \in \{6,7,23,24,29,30,33,74\}, \\
0 & \textrm{ otherwise, } 
\end{array}
\right.
$$
where $A \ = \{1,2,5,11,12,13,16,17,19,20,21,25,37,45,55,65,68,70,71,82,84,
\newline 89,98,112,123,130,135,
151,189,198,203,220,245,267,345,355,359,381,401,
\newline 443,464,514,561,583,610,612,624\}.$ As said above, under this particular morphism, no Pythagorean triple in the interval $[1,4632]$ is monochromatic. For the record, there are 735 primitive Pythagorean triples in that interval.

\bigskip
Propositions~\ref{morphism m=2} and~~\ref{morphism m=3} give rise to an interesting problem which might be more manageable than the general regularity problem of the Pythagoras equation.
\begin{problem} Is it true that, for any $m \ge 4$ and any standard morphism $f \colon \N_+ \to \Z/m\Z$, monochromatic Pythagorean triples are unavoidable? And if yes, how does the unavoidability threshold behave as a function of $m$?
\end{problem}
Recall that for $m=2$ and $3$, the corresponding unavoidability threshold turns out to be $533$ and $4633$, respectively.

\section{Morphic colorings}
We now turn to $\PN_0$-partial morphisms in the sense of Section~\ref{sect3}.
\begin{definition} A \emph{morphic coloring} of $\N_+$ is a $\PN_0$-morphism where $\PN_0$ is a finite subset of $\PN$. 
\end{definition}
In the sequel, we shall mostly consider morphic colorings of $\N_+$ in this sense.
Let us start first with an easy remark concerning primitive Pythagorean triples.

\begin{remark} There exist partial morphisms $f \colon \N_+ \to \Z/2\Z$ under which no \emph{primitive} Pythagorean triple is monochromatic.
\end{remark}
One obvious example is the 2-coloring given by $f(n)= 1$ if $n$ is even and 0 otherwise. Since any primitive Pythagorean triple $(a,b,c)$ contains exactly one even number, it is not monochromatic under $f$. 

Two similar examples arise by mapping multiples of 3 to color 1, or else multiples of 5 to color 1, and the rest to color 0, respectively. Indeed, any primitive Pythagorean triple $(a,b,c)$ contains at least one multiple of 3, and one multiple of 5 as well; this easily follows from the fact that the only nonzero square mod 3 is 1, and the only nonzero squares mod 5 are $\pm 1$. But since $a,b,c$ are assumed to be coprime, and since $a^2+b^2=c^2$, they cannot be \textit{all three} mapped to 1, or to 0, by these two 2-colorings.

These three 2-colorings are $\emptyset$-morphisms, as they satisfy $f(xy)=f(x)+f(y)$ for all coprime positive integers $x,y$. More precisely, they are characterized by the values
$$
f(p_0^\nu)= 1, \ f(p^\nu)=0
$$
for all $\nu \ge 1$ and all primes $p \not = p_0$, where $p_0 = 2$, 3 or 5, respectively.

\medskip

Perhaps surprisingly, for any $\emptyset$-morphism $f \colon \N_+ \to \Z/2\Z$, it turns out again that monochromatic Pythagorean triples are unavoidable. Here are even stronger results, obtained by exhaustive computer search with an algorithm briefly described below.

\begin{proposition} For any $\PN_0$-morphism $f \colon \N_+ \to \Z/2\Z$, where $\PN_0$ is one of the sets $\{2,3,5\}$, $\{2,3,5,7\}$ and $\{2,3,5,7,11\}$, monochromatic Pythagorean triples are unavoidable. More precisely, they are unavoidable in the integer interval $[1,N]$, where
$$
N \ = \ 
\left\{
\begin{array}{rcl}
533 & \textrm{if} & \PN_0=\{2,3,5\}, \\
565 & \textrm{if} & \PN_0=\{2,3,5,7\}, \\
696 & \textrm{if} & \PN_0=\{2,3,5,7,11\}.
\end{array}
\right.
$$
Moreover, in each case, the given $N$ is minimal with respect to that property.
\end{proposition}

\subsection{The algorithm}

Here is a brief description of the algorithm used. It consists of a recursive, backtracking 
search, that tries to color all elements in Pythagorean triples within a given integer interval $[1,M]$ without creating monochromatic such triples. 

Let us consider a fixed subset $\PN_0 \subset \PN$. The set of variables is then the set $S(\PN_0)$ of positive integers defined in Section~\ref{sect3}, namely
$$
S(\PN_0) = \big\{n_0 \in \N_+ \mid \textrm{supp}(n_0) \subseteq \PN_0\big\} \, \bigsqcup \, \big\{p^\nu \mid p \in \PN \setminus \PN_0, \nu \in \N_+\big\}.
$$
For $n \in \N_+$, let us denote by
$\f(n)=\{q_1, \ldots, q_k\}$ the unique subset of $S(\PN_0)$ such that
$$n = \prod_{i=1}^{k} q_i$$
and where each $q_i \in S(\PN_0)$ is \emph{maximal}, in the sense that no proper multiple of $q_i$ dividing $n$ belongs to $S(\PN_0)$. Thus, for a given $\PN_0$-morphism $f \colon \N_+ \to \Z/2\Z$, and for $n \in \N_+$, we have
$$
f(n) = \sum_{q \in \f(n)} f(q).
$$
For instance, if $\PN_0=\{2,3,5\}$ and $n=64680=2^3\cdot 3\cdot 5 \cdot 7^2 \cdot 11$, the maximal $S(\PN_0)$-factors of $n$ are $120=2^3\cdot 3\cdot 5$, $49=7^2$ and $11$. Thus $\f(n)=\{120,49,11\}$, and $f(n)=f(120)+f(49)+f(11)$ for any $\PN_0$-morphism $f$ as above.

\medskip

The algorithm will try to assign a suitable color in $\Z/2\Z$ to each variable, but the order in which this is done is important and may strongly affect the running time. To define a proper assignment order, we introduce the following notation. Given a positive integer $M$, let $\T_M$ denote the set of all Pythagorean triples contained in the integer interval $[1,M]$. Then, for $q \in S(\PN_0)$ and ${\{a,b,c\} \in \T_M}$, we define
$$
\delta_{q}^{\{a,b,c\}} \ = \
\left\{
\begin{array}{rl}
  1 & \textrm{ if } q \in \f(a) \cup \f(b) \cup \f(c), \\
  0 & \textrm{ otherwise. } 
\end{array}
\right.
$$
The \emph{weight} of the variable $q$ is now defined as
$$w(q) = \sum_{t \in \T_M}\delta_{q}^{t},$$
i.e. the number of Pythagorean triples in $[1,M]$ where $q$ appears in the 
decomposition $\f(n)$ of one of the triple elements $n$. The variables are then ordered by decreasing weight, and the algorithm assigns a value in 
$\Z/2\Z$ to the variables in that order. Thus, variables constrained by the greatest number of triples in which they are involved as a maximal $S(\PN_0)$-factor are tested first.

Once a variable is assigned, the algorithm performs forward arc
checking~\cite{bessiere2006constraint}; that is, it computes, if possible, the color of all 
triple elements following the current
partial morphism. This color computation of an element $m$ is possible
if all variables in $\f(m)$ are already colored as explained in
Section~\ref{sect3}.  If the coloration of any member creates a
monochromatic triple, then the algorithm tries the other color for the
variable if any left, or backtracks.

\subsection{The function $N(\PN_0)$}

The above results prompt us to introduce the following function.

\begin{definition} For any subset $\PN_0 \subseteq \PN$, we denote by $N(\PN_0)$ the largest integer $N$ if any, or $\infty$ otherwise, such that there exists a 2-coloring of the integer interval $[1,N]$ by a $\PN_0$-morphism which avoids monochromatic Pythagorean triples in that interval.
\end{definition}

The above results may thus be expressed as follows:
\begin{eqnarray*}
N(\{2,3,5\}) & = &  532, \\
N(\{2,3,5,7\}) & = &  564, \\
N(\{2,3,5,7,11\}) & = &  695.
\end{eqnarray*}

\smallskip

Moreover, whether the equation $X^2+Y^2=Z^2$ is 2-regular or not is equivalent to whether $N(\PN)$ is finite or infinite, respectively. This follows from a standard compactness argument.

\subsection{On small sets of primes up to 100}\label{up to 100}
Finally, we consider $\PN_0$-morphisms $f \colon \N_+ \to \Z/2\Z$ where $\PN_0$ ranges through all sets of prime numbers in $[2,100]$ of cardinality 3, 4 and 5. Note that $[2,100]$ contains 25 prime numbers. Needless to say, a considerable amount of code optimization and computer time were needed in order to establish the findings below. 

The interest of the results below is that the unavoidability thresholds of monochromatic Pythagorean triples remain much smaller than $7825$, the general threshold of Theorem~\ref{heule}, and also that their numerical behavior turn out to be quite subtle.

\subsubsection{The case $|\PN_0|=3$}

There are $\displaystyle \binom{25}{3} = 2300$ sets of three distinct prime numbers smaller than 100. 

\begin{proposition}\label{triples} Among the $2300$ subsets $\PN_0 \subset \PN \cap [2,100]$ of cardinality 3, one has $N(\PN_0)=532$ for all but 29 of them. These 29 exceptions are as follows:
$$
N(\PN_0) \ = \ 
\left\{
\begin{array}{rl}
544 & \textrm{ if } \ \PN_0 \in \big\{\{2,3,13\},\{3,5,17\},\{5,13,41\}\big\}, \\
564 & \textrm{ if } \ \PN_0 \in \big\{\{2,3,7\}\big\} \cup \{\{7,11,a\} \mid a \in \PN \cap [2,100]\setminus \{7,11\}\big\}, \\
628 & \textrm{ if } \ \PN_0 \in \big\{\{2,3,19\},\{3,13,19\}\big\}.
\end{array}
\right.
$$
\end{proposition}

\subsubsection{The case $|\PN_0|=4$}
There are $\displaystyle \binom{25}{4} = 12650$ sets of four distinct prime numbers smaller than 100.
\begin{proposition} For all subsets $\PN_0 \subset \PN \cap [2,100]$ of cardinality 4, one has
$$
532 \ \le \ N(\PN_0)  \ \le \ 784
$$
and, more precisely, 
\begin{eqnarray*}
N(\PN_0) & \in & \{532, 543, 544, 547, 564, 577, 594, 614, 624, 628, \\ 
& & 649, 656, 662, 666, 679, 688, 696, 739, 778, 784\}.
\end{eqnarray*} 
\end{proposition}
\begin{remark} The value $784$ above is attained only once, by $\PN_0=\{3,5,19,23\}$. The same holds for the next few largest values, including $778=N(\{3,13,19,23\})$ and $739=N(\{2,3,19,23\})$. All three cases involve the subset $\{3,19,23\}$. Compare with the most performant triples found in Proposition~\ref{triples}, namely $\{2,3,19\}$ and $\{3,13,19\}$, both containing $\{3,19\}$.
\end{remark}
\begin{remark}
The behavior of $N(\PN_0)$ may be quite subtle. For instance, one has $
N(\PN_0) = 564$ for all quadruples $\PN_0$ in $\PN \cap [2,100]$ containing the pair $\{7,11\}$, with one single exception given by $N(\{7,11,13,17\}) = 624$. 
\end{remark}

\subsubsection{The case $|\PN_0|=5$}
There are $\displaystyle \binom{25}{5} = 53130$ sets of five distinct prime numbers smaller than 100.

\begin{proposition} For every subset $\PN_0 \subset \PN \cap [2,100]$ of cardinality 5, we have
$$
532 \ \le \ N(\PN_0)  \ \le \ 900.
$$
Moreover, the only such subsets $\PN_0$ attaining the maximum $N(\PN_0) = 900$ are
$$
\{2,3,7,19,23\}, \ \{2,3,17,19,23\}. 
$$
\end{proposition}

\subsection{Summary}
The findings of section~\ref{up to 100} may be summarized as follows. 
\begin{proposition} For every subset $\PN_0 \subset \PN \cap [2,100]$ of cardinality at most 5, and for every 2-coloring
$$
f \colon [1,901] \to \Z/2\Z
$$
by a $\PN_0$-morphism, monochromatic Pythagorean triples are unavoidable.
\end{proposition}

\medskip

As for the $\displaystyle \binom{25}{6} = 177100$ subsets of cardinality 6 of $\PN \cap [2,100]$, an exhaustive search cannot currently be completed in a reasonable amount of time. A first partial search has yielded 900 again as the highest value of $N(\PN_0)$ found so far, achieved by the following subsets:
$$
\begin{array}{c}
\{2, 3, 5, 7, 19, 23\}, \,
\{2, 3, 5, 11, 19, 23\}, \,
\{2, 3, 5, 19, 23, 41\}, \, \\
\{2, 3, 5, 17, 19, 23\}, \,
\{2, 3, 5, 19, 23, 47\}, \,
\{2, 3, 5, 19, 23, 53\}.
\end{array}
$$
Interestingly, looking more closely at the results of section~\ref{up to 100}, one notes that the subsets $\PN_0 \subset \PN \cap [2,100]$ maximizing the function $N$ all contain $\{3,19\}$ if $|\PN_0|=3$, or $\{3,19,23\}$ if $|\PN_0|=4$, or $\{2,3,19,23\}$ if $|\PN_0|=5$. While the case $|\PN_0|=6$ is largely incomplete, an analogous statement might well hold.

\bigskip
\noindent

We thank the Editors of this journal and the anonymous referee for viewing favorably this paper in light of the general $2$-regularity result of \cite{HKM:15}, and the referee again for his/her highly valuable comments. We also thank Gyan Prakash for interesting information related to this work and for reference \cite{FH:14}.

\bigskip
\bigskip

\noindent
{\small
\textbf{Authors addresses:}

\bigskip

\noindent
$\bullet$ Shalom Eliahou, Jean Fromentin\textsuperscript{a,b}:

\noindent
\textsuperscript{a}Univ. Littoral C\^ote d'Opale, EA 2597 - LMPA - Laboratoire de Math\'ematiques Pures et Appliqu\'ees Joseph Liouville, F-62228 Calais, France\\
\textsuperscript{b}CNRS, FR 2956, France\\
e-mail: \{eliahou, fromentin\}@lmpa.univ-littoral.fr, \{eliahou, fromentin\}@math.cnrs.fr

\bigskip

\noindent
$\bullet$ Virginie Marion-Poty, Denis Robilliard\textsuperscript{a}:

\noindent
\textsuperscript{a}Univ. Littoral C\^ote d'Opale, EA 4491 - LISIC - Laboratoire d'Informatique Signal et Image de la C\^ote d'Opale, F-62228 Calais, France\\
e-mail: \{poty, robilliard\}@lisic.univ-littoral.fr

}

\end{document}